\documentclass[11pt,oneside,english]{amsart}
\usepackage[T1]{fontenc}
\usepackage[latin9]{inputenc}
\usepackage{geometry}
\usepackage{verbatim}
\usepackage{float}
\usepackage{amstext}
\usepackage{amsthm}
\usepackage{amssymb}
\usepackage{mathdots}
\usepackage{graphicx}

\makeatletter

\floatstyle{ruled}
\newfloat{algorithm}{tbp}{loa}
\providecommand{\algorithmname}{Algorithm}
\floatname{algorithm}{\protect\algorithmname}

\numberwithin{equation}{section}
\numberwithin{figure}{section}
\theoremstyle{plain}
\newtheorem{thm}{\protect\theoremname}
  \theoremstyle{plain}
  \newtheorem{cor}[thm]{\protect\corollaryname}
  \theoremstyle{plain}
  \newtheorem{lem}[thm]{\protect\lemmaname}

\usepackage{tikz}
\usepackage{algpseudocode}
\usepackage{multicol}

\makeatother

\usepackage{babel}
  \providecommand{\corollaryname}{Corollary}
  \providecommand{\lemmaname}{Lemma}
\providecommand{\theoremname}{Theorem}

\begin{document}

\title{Recovering the Structure of Random Linear Graphs}

\author{Israel Rocha}

\author{Jeannette Janssen}

\author{Nauzer Kalyaniwalla}

\thanks{israel.rocha@dal.ca, jeannette.janssen@dal.ca. Department of Mathematics
and Statistics. nauzerk@cs.dal.ca. Faculty of Computer Science. Dalhousie
University, Halifax, Canada.}
\begin{abstract}
In a random linear graph, vertices are points on a line, and pairs
of vertices are connected, independently, with a link probability
that decreases with distance. We study the problem of reconstructing
the linear embedding from the graph, by recovering the natural order
in which the vertices are placed. We propose an approach based on
the spectrum of the graph, using recent results on random matrices.
We demonstrate our method on a particular type of random linear graph.
We recover the order and give tight bounds on the number of misplaced
vertices, and on the amount of drift from their natural positions.
\end{abstract}

\maketitle

\section{Introduction}

Spatial networks are graphs with vertices located in a space equipped
with a certain metric. In a random spatial network two vertices are
connected if their distance is in a given range. Random spatial networks
have been used in the study of DNA reconstruction, very large scale
integration (VLSI) problems, modelling wireless adhoc networks, matrix
bandwidth minimization, etc. The book \cite{Penrose} contains a rich
source for the mathematics behind these structures. 

One particular instance of these graph models are the so called one-dimensional
geometric graphs, where the vertices are points in $\mathbb{R}$,
connected with some probability if they are close. This is the seriation
problem \cite{Atkins}. Seriation is an exploratory data analysis
technique to reorder vertices into a sequence along a line so that
it reveals its pattern. This is an old problem tracing back to the
19-th century and it arises in many different areas of science. For
a historical overview on this subject we recommend \cite{Liiv}. Also,
this simple graph model was successfully used to predict key structural
properties of complex food webs \cite{Williams} and protein interaction
networks \cite{Pires}.

In this paper we work with this kind of graph model, which we call
random linear graphs. We are concerned with large amounts of vertices
$n$, where they are connected with some probability if their distance
is at most $n/2$. We show how to successfully retrieve the linear
order of vertices that are randomly distributed. Besides, we show
that this order is correct with high probability, and we quantify
how many vertices are misplaced. That is the first time one fully
recovers the structure of random linear graphs using this method and
it serves as a proof of concept. In a forthcoming work we will deal
with different models using the same technique.

Closely related to our work, in \cite{Fogel} the authors considered
the problem of ranking a set of $n$ items given pairwise comparisons
between these items. They showed that the ranking problem is directly
related to the seriation problem. That allows them to provide a spectral
ranking algorithm with provable recovery and robustness guarantees.

Another model that is related to random linear graphs is the Stochastic
Block Model. There, the set of vertices is partitioned in disjoint
sets of communities $C_{1},\ldots,C_{k}$, and vertices $u\in C_{i}$
and $v\in C_{j}$ are connected by an edge with probability $p_{ij}$.
If one consider the communities to be formed by one single vertex,
and fix $p_{ij}=p$ if the distance between $v_{i}$ and $v_{j}$
is at most $n/2$, then we obtain the random linear graph model. The
difference from our case is that the Stochastic Block Model makes
vertices inside the same community indistinguishable. In our case,
each vertex has its own place - corresponding to a community in the
Stochastic Block Model. That makes the problem of recovering the full
structure of the linear model harder than recovering the structure
in the Stochastic Block Model.

The main tools we use to recover the structure of such graphs are
the eigenvectors of its adjacency matrix. Eigenvectors were successfully
applied in the past for relevant problems such as recovering a partition,
clique, coloring, bipartition, etc., which are naturally present,
but hidden in a random graph \cite{Alon,AlonSudakov,Bopana,Bui,Condon,Dekel}.
An approach that is closely related to the technique we use, is the
spectral partitioning in the Stochastic Block Model. In one of the
first papers to apply this idea \cite{McSherry}, McSherry shows how
to solve the hidden partition problem for certain parameters. His
investigation has since been improved and we can find recent developments
in \cite{ChinVu,Vu-1} for instance. 

McSherry's technique relies on the low rank property of the Stochastic
Block Model matrix. For this model, the rank of the matrix is equal
to the number of blocks. In order to distinguish vertices that belong
to different blocks, this technique uses a number of eigenvectors
equal to the rank of the matrix. Thus, for a model with few blocks,
we only need a few eigenvectors. The intuition is that the information
necessary to distinguish vertices inside different blocks is encoded
in the top singular vectors of the model matrix. That is because they
provide a good low rank approximation for the model matrix. 

For the random linear graph model something different happens. Its
matrix has full rank, thus the idea of using a few eigenvectors to
provide a good approximation seems hopeless. There is where our method
begins. Instead of approximating the model matrix by singular vectors,
we identify which eigenvectors encode the structure of the graph itself.
We show that the structure of the graph in question is encoded inside
only one eigenvector of its adjacency matrix. Therefore, that eigenvector
is the important object to understand. In the next section we formalize
the problem. We show how to recover the structure of random linear
graphs.

%The first step of our approach is to find eigenvectors that encode
%the linear structure of the graph model. With that, we bound the distance
%between the eigenvectors of the random graph and the graph model.
%Then results from random matrix theory
%are not enough to provide good bounds. The quality of our bounds relies
%on more precise informations about the spectrum of the graph model.
%We need to bound the gap between eigenvalues to be able to provide
%good results. Besides, we want to quantify the error in the ordering
%provided by our approach. To see how the error in the vector approximation
%translates into error in the ordering, we need closed expressions
%for the eigenvectors of model matrix. 
Turns out the matrices that arise from the random linear graph model are certain Toeplitz matrices
which have unknown spectrum so far. Besides the fact that Toeplitz
matrices are well studied in the literature for many years, there
is a lack of closed expressions for its spectrum in the general form.
To convince the reader about its difficulty we provide the references
\cite{boettcher,Gray}.

The problem we address can be formulated as follows: given a graph that is an instance of a linear random graph process, extract the ordering of the positions of the vertices on the line. To solve this problem, we rely on two tools. First, we use spectral theory to show that the correct ordering can be extracted from the eigenvectors of the model matrix. Then, we use results from random matrix theory to show that the eigenvectors of the model matrix and the adjacency matrix of the graph are similar. Finally, we derive how this similarity implies that the ordering from the eigenvectors of the adjacency matrix closely matches the true ordering. 

The main challenge of this approach is that it requires detailed knowledge of the spectrum of the model matrix. First, the bound on the similarity between eigenvectors of model matrix and adjacency matrix requires lower bounds on the gaps between eigenvalues. More importantly, to derive how the bound on the difference of the eigenvectors of model matrix and adjacency matrix translates into bounds on the errors in the ordering, we need to have precise knowledge about the eigenvector. Specifically, let $v$ be an eigenvector of the model matrix, and assume that $v$ reveals the true ordering, i.e.~if the vertices are ordered so that the corresponding components of $v$ are increasing, then this ordering corresponds to the linear embedding. Let $\hat{v}$ be the corresponding eigenvector of the adjacency matrix, and assume that $||v-\hat{v}||$ is small. Without further knowledge, we cannot conclude that the ordering obtained from
$\hat{v}$ is close to the true ordering. In order to conclude that a small perturbation of $v$ cannot lead to large changes in the ordering, the components of $v$ must not only be increasing, but successive components of $v$ must be relatively far apart. In other words, $v$ must be ``steeply increasing."

%In order to illustrate our approach, we have chosen one particular type of linear graph model.
%Namely, this is a graph model where vertices are placed at equal intervals in the line segment $[0,1]$, and they are made adjacent with a given probability $p$ if they are less than $0.5$ apart. We show that, for this model, the second eigenvector of the adjacency matrix gives an ordering that is arbitrarily close (in terms of Kendall's tau coefficient) to the true ordering. We also give precise bounds on the number of vertices that are misplaced.

\section{Problem statement and main results}

The random linear graph in the class we consider will be denoted as $G=(V,E)$, with
set of vertices $\left\{ v_{1},v_{2}\ldots,v_{n}\right\} $, where
we put an edge between each pair of vertices $v_{i}$ and $v_{j}$
with probability $p$ if $\left|i-j\right|\leq n/2$. The {\sl model matrix} (matrix of edge probabilities) is a banded matrix that describes a graph with
a clear linear structure. 

\begin{figure}[H]
\begin{tabular}{ccc}
\begin{tikzpicture}[scale=0.05] \draw[fill] (0,0) circle [radius=1]; \draw[fill] (10,0) circle [radius=1]; \draw[fill] (20,0) circle [radius=1]; \draw[fill] (30,0) circle [radius=1]; \draw[fill] (40,0) circle [radius=1]; \draw[fill] (50,0) circle [radius=1]; \draw[fill] (60,0) circle [radius=1]; \draw  (20,0) arc [radius=11.4, start angle=25, end angle= 155]; \draw  (30,0) arc [radius=11.4, start angle=25, end angle= 155]; \draw  (40,0) arc [radius=11.4, start angle=25, end angle= 155]; \draw  (50,0) arc [radius=11.4, start angle=25, end angle= 155]; \draw  (60,0) arc [radius=11.4, start angle=25, end angle= 155]; \draw  (30,0) arc [radius=15, start angle=5, end angle= 175]; \draw  (40,0) arc [radius=15, start angle=5, end angle= 175]; \draw  (50,0) arc [radius=15, start angle=5, end angle= 175]; \draw  (60,0) arc [radius=15, start angle=5, end angle= 175]; \draw (0,0) -- (60,0); \end{tikzpicture}   &   {\tiny $M=p\left[\begin{array}{ccccccc} 0&\color{blue}{1}&\color{blue}{1}&\color{blue}{1}&0&0&0\\ \color{blue}{1}&0&\color{blue}{1}&\color{blue}{1}&\color{blue}{1}&0&0\\ \color{blue}{1}&\color{blue}{1}&0&\color{blue}{1}&\color{blue}{1}&\color{blue}{1}&0\\ \color{blue}{1}&\color{blue}{1}&\color{blue}{1}&0&\color{blue}{1}&\color{blue}{1}&\color{blue}{1}\\ 0&\color{blue}{1}&\color{blue}{1}&\color{blue}{1}&0&\color{blue}{1}&\color{blue}{1}\\ 0&0&\color{blue}{1}&\color{blue}{1}&\color{blue}{1}&0&\color{blue}{1}\\ 0&0&0&\color{blue}{1}&\color{blue}{1}&\color{blue}{1}&0\\ \end{array}\right]$}\\
\end{tabular} 

\caption{Model graph and its model matrix}
\end{figure}
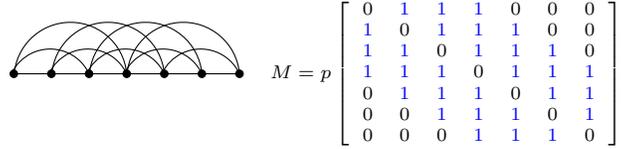

Let $M$ be this model matrix with entries $m_{ij}=p$ if $\left|i-j\right|\leq n/2$
and $0$ otherwise. Furthermore, the adjacency matrix $\widehat{M}$
of the random linear graph is a random matrix whose entries are independent
Bernoulli variables, where $\mathbb{P}\left(\widehat{m}{}_{ij}=1\right)=m_{ij}$. 

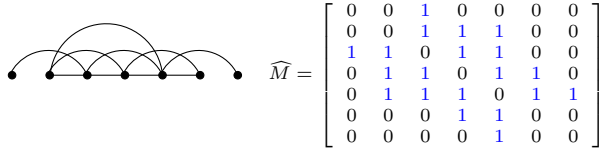
\begin{figure}[H]
\begin{tabular}{ccc}
\begin{tikzpicture}[scale=0.05] \draw[fill] (0,0) circle [radius=1]; \draw[fill] (10,0) circle [radius=1]; \draw[fill] (20,0) circle [radius=1]; \draw[fill] (30,0) circle [radius=1]; \draw[fill] (40,0) circle [radius=1]; \draw[fill] (50,0) circle [radius=1]; \draw[fill] (60,0) circle [radius=1]; \draw  (20,0) arc [radius=11.4, start angle=25, end angle= 155]; \draw  (30,0) arc [radius=11.4, start angle=25, end angle= 155]; \draw  (40,0) arc [radius=11.4, start angle=25, end angle= 155]; \draw  (50,0) arc [radius=11.4, start angle=25, end angle= 155]; \draw  (60,0) arc [radius=11.4, start angle=25, end angle= 155]; \draw  (40,0) arc [radius=15, start angle=5, end angle= 175]; \draw (10,0) -- (50,0); \end{tikzpicture}   &   {\tiny $\widehat{M}=\left[\begin{array}{ccccccc} 0&0&\color{blue}{1}&0&0&0&0\\ 0&0&\color{blue}{1}&\color{blue}{1}&\color{blue}{1}&0&0\\ \color{blue}{1}&\color{blue}{1}&0&\color{blue}{1}&\color{blue}{1}&0&0\\ 0&\color{blue}{1}&\color{blue}{1}&0&\color{blue}{1}&\color{blue}{1}&0\\ 0&\color{blue}{1}&\color{blue}{1}&\color{blue}{1}&0&\color{blue}{1}&\color{blue}{1}\\ 0&0&0&\color{blue}{1}&\color{blue}{1}&0&0\\ 0&0&0&0&\color{blue}{1}&0&0\\ \end{array}\right]$}\\
\end{tabular} 

\caption{Random linear graph and its random matrix}
\end{figure}

\begin{figure}[H]
\includegraphics[scale=0.3]{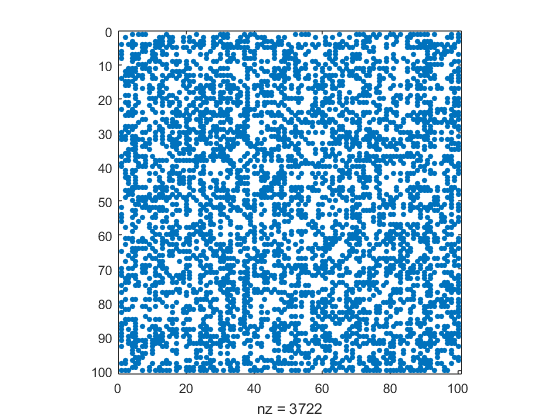}\includegraphics[scale=0.3]{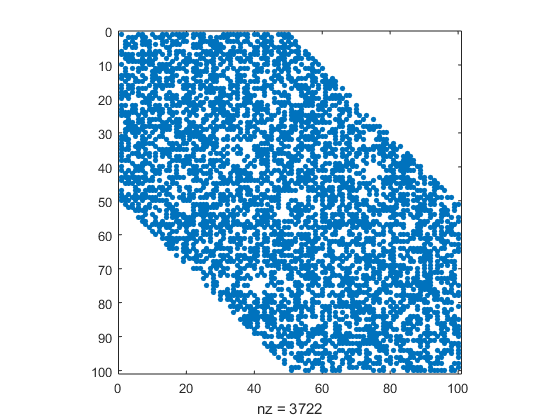}

\caption{Random matrix before and after an appropriate permutation.}
\label{fig:rearrangedMatrix}
\end{figure}

For the random linear graph in the last figure, the order in which
the vertices appear makes its linear structure clear, and this is
reflected in the band structure of its adjacency matrix. For large
matrices with unknown ordering it might be a challenge to reveal its
correct linear structure, as Figure \ref{fig:rearrangedMatrix} shows.

Now we can formulate the problem of reconstructing the linear embedding
of the random graph as follows. Given a graph produced by a linear
random graph model, find its linear embedding, or equivalently, retrieve
the linear order of the vertices. 

The model graph is a unit interval graph. It is well-known that the linear order can be retrieved in polynomial time, see for example \cite{Booth,Corneil1,Corneil2}.
The random graph we consider can be seen as a random subgraph of this graph. Such a graph has a clear linear structure, but the linear ordering cannot be retrieved with the same methods. 

In the particular graph we consider, the order of the vertices of the model graph can also be retrieved by looking at the degree of the vertices: the end vertices have degree $n/2$, and the degrees increase near the middle, with a maximum of $n-1$ for the center vertices. 

The degree of a vertex in the random graph is the result of a large number of idenpendent trials with identical probability. 
Thus, one could try to estimate the position of the vertices in the random matrix by considering the degree (vertices from both ends can be distinguished by the number of common neighbours). With high probability, the degree of vertex $v_i$, is 
$$(n/2+i-1)p + \mathcal{O}(\sqrt{n}\log n), \text{ for }  i=1,\ldots, n/2,\text{ and }$$
$$(3n/2-i+1)p + O(\sqrt{n}\log n), \text{ for } i=n/2+1,\ldots,n.$$
Thus, most vertices will drift at most $O(\sqrt{n} \log n)$ positions from their correct position, and consequently the number of inverted pairs is $o(n^2)$. However, one can expect a non-neglible portion of all pairs to have a drift of $\sqrt{n}$. In comparison, our method gives the  correct ordering for all vertices, except a small fraction near the ends. Moreover, in our approach, we can provide precise bounds on the number of pairs out of order and on the distance of vertices from their correct positions. As we will see, the results exhibits a trade off between how far vertices are from their original position and how far they are from the end, and
the total number of such incorrectly placed vertices. 
%Notice that the first and last quarter of vertices do not have common neighbors. In view of that, we can determine the position of the first and last quarter of vertices up to an error of $\mathcal{O}(\sqrt{n})$ using its degrees. However, vertices in the middle have a significant overlap in their neighborhood, thus the degrees alone are not enough to determine if their positions are in the first or second half. That makes it a hard approach to distinguish in which position vertices belong to.

It is also possible to use a Stochastic Block Model to approximate the model matrix $M$ with blocks of size $\sqrt{n}$, for example. But that again gives only an approximated ordering, where vertices that belong to the same block will be positioned anywhere within it.  

\begin{figure}
\begin{tabular}{cccc} $p=0.1$ & $p=0.5$ & $p=0.9$\\ \includegraphics[width=0.3\textwidth]{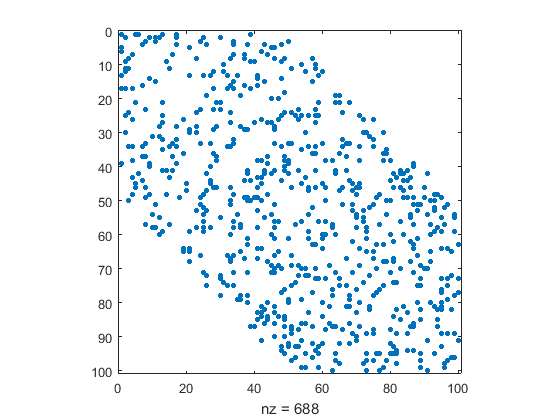}& \includegraphics[width=0.3\textwidth]{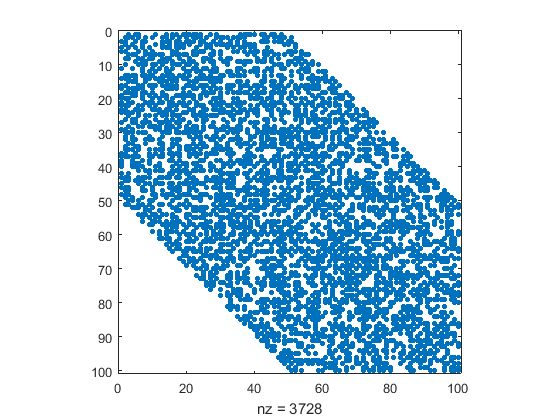}& \includegraphics[width=0.3\textwidth]{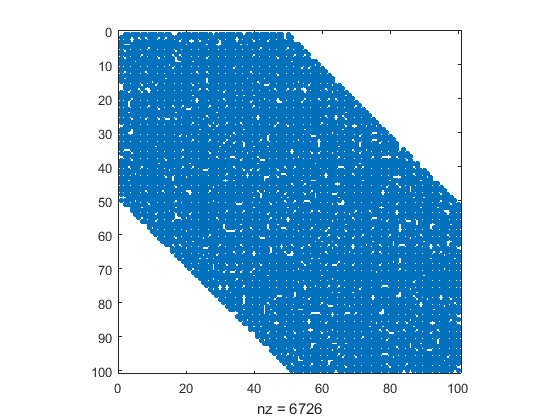}\\ \end{tabular}
\begin{tabular}{ccc}  \includegraphics[width=0.3\textwidth]{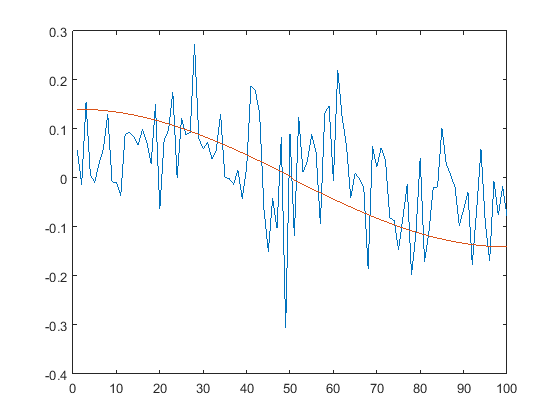} & \includegraphics[width=0.3\textwidth]{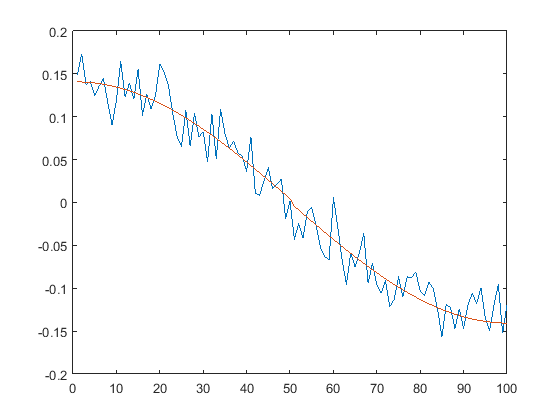} & \includegraphics[width=0.3\textwidth]{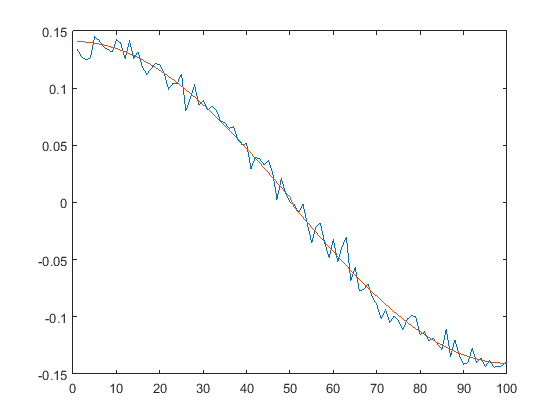} \\ \end{tabular}

\caption{}
\label{fig:matrices}
\end{figure}

Before presenting an algorithm that
retrieves the correct order we look at the eigenvectors
of the matrices in question. We plot a specific eigenvector of the adjacency
matrix by plotting the component values as a  function of the indices. Figure \ref{fig:matrices} shows,
for different probabilities $p$, the random matrix in the correct
order and the eigenvector of the second largest eigenvalue.

The smooth curve in the graph of Figure \ref{fig:matrices} represents
the eigenvector of the deterministic model matrix in the correct order.
We will prove that the entries of this eigenvector are always monotonic
for different matrix dimensions. Since its entries are monotonic, one
could use it to find the correct embedding even when the matrix is
given with the wrong order. 

In Figure \ref{fig:matrices}, the other
line represents the entries of the eigenvector of the random matrix.
Notice that the spiked lines exhibit the same decreasing trend as
the straight lines. We will show that this is indeed true with high
probability. That is, the corresponding eigenvector of the random matrix gives
an order that is close to the correct order with high probability.
Our first theorem shows that the second eigenvectors of random matrix and model 
matrix are close in norm.
\begin{thm}
\label{thm:upperBound}Let $M$ be the linear graph model matrix with
constant probability $p$ and variance $\sigma^{2}$, and let $\widehat{M}$
the random matrix following the model matrix. Let $x$ be a unitary
eigenvector for $\lambda_{2}(M)$ and $\widehat{x}$ be an unitary
eigenvector for $\lambda_{2}(\widehat{M})$. There is a constant $C_{0}>0$
such that 
\[
\left\Vert x-\widehat{x}\right\Vert \leq C_{0}n^{-1/2},
\]
with probability at least $1-n^{-3}$.
\end{thm}
%Notice that the higher the probability, the closer the eigenvectors
%are, which portray the behavior displayed in Figure \ref{fig:matrices}.
%jj: I removed this sentence, because the theorem does not mention $p$, so this
%observation, while true, does not follow from the theorem as stated.
Theorem \ref{thm:upperBound} suggests the following algorithm
to recover the order of a random linear graph.

\begin{algorithm}[H]
\begin{algorithmic} \Require{Random matrix $\widehat{M}$} \\Compute $\widehat{x}$, the eigenvector for $\lambda_{2}(\widehat{M})$ \\
Permute $\widehat{M}$ according to the order of the entries of $\widehat{x}$  \\
\Return{$\widehat{M}$} \end{algorithmic}

\caption{Order the adjacency matrix of the random linear graph }
\label{alg:algo1}
\end{algorithm}

One important question that we are concerned with is the quality
of the order provided by Algorithm \ref{alg:algo1}. In particular, we
want to know how many of such vertices exist and how far they are
from their correct positions (see Section \ref{sec:qualitative}).
We show that there is a trade off between how far vertices can be
misplaced and the total amount of misplaced vertex pairs. In Section \ref{sec:qualitative}
we also calculate the rank correlation coefficients for the order
provided by Algorithm \ref{alg:algo1}. 

A rank correlation coefficient measures the degree of similarity between
two lists, and can be used to assess the significance of the relation
between them. One can see the rank of one list as a permutation of
the rank of the other. Statisticians have used a number of different
measures of closeness for permutations. Some popular rank correlation
statistics are Kendall's $\tau$, Kendall distance, and Spearman's
footrule. There are several other metrics, and for different situations
some metrics are preferable. For a deeper discussion on metrics on
permutations we recommend \cite{Diaconis}. 

In this paper, we will use the metrics defined as follows. The permutation $\sigma$
 derived from an $n$-dimensional vector $y$ is the permutation obtained from ordering the components of $y$ in decreasing order. Specifically, we have that
$y_{\sigma(1)}\geq y_{\sigma(2)}\geq\ldots\geq y_{\sigma(n)}$.
When convenient we indicate it by $\sigma_{y}$. 

To count the total number of inversions in $\sigma$ we use
\[
D(\sigma)=\sum_{i<j}\textbf{1}_{\sigma(i)>\sigma(j)}\text{ (Kendall Distance)},
\]
The distance that an element $i$ moved due to the permutation (the {\em displacement or drift})
is $\left|i-\sigma(i)\right|$. To count the total drift, we use
\[
F(\sigma)=\sum_{i=1}^{n}|i-\sigma(i)|\text{ (Spearman's footrule)},
\]
which is the total displacement  of all elements. Finally, Kendall's $\tau$
rank correlation coefficient is defined as 

\[
\tau=\frac{(\#concordant\ pairs-\#discordant\ pairs)}{\binom{n}{2}}=1-\frac{D(\sigma)}{n(n-1)}.
\]

Even though metrics can be quite different from each other, in \cite{DiaconisGraham}
Diaconis and Graham proved that Kendall distance and Spearman's Footrule
are equivalent measures in the following sense:

\[
D(\sigma)\leq F(\sigma)\leq2D(\sigma).
\]
Thus, if the lists to be compared are large enough, these measures
are of the same order. That is exactly our situation, where we want
to quantify the number of vertices in the wrong order as $n$ tends
to infinity.

Before providing metrics for the order given by the eigenvector of
the random matrix, we need a better understanding of the displacement
of vertices in that order. More precisely, we need to quantify how
far vertices can be from the correct position and how many of them
exists. In fact, we will see that this order is very accurate as the
number of vertices grows.

To this end, we find an expression for the minimum distance between
eigenvector $y$ of the random matrix and the known eigenvector $x$
of the model. We show that if there are too many incorrectly placed
vertices, this lower bound is greater than the upper bound of Theorem
\ref{thm:upperBound}. That means $y$ is too far from $x$, which
happens with small probability. Therefore, the number of vertices
in the wrong position can be quantified with high probability.

The next theorem reveals how the permutation of the eigenvector of
the random matrix $\widehat{M}$ compares to the eigenvector of $M$.
It says that large portions of vertices that are in the wrong order
cannot be too far away from their correct positions, and entries that
are too far represent a small group of vertices.

First we define a refined version of the Kendall distance. This version
counts inverted pairs that appear after position $r$, and whose indices are at least $k$ positions apart.
First note that, for a vector $y$, we can rewrite $D(\sigma_{y})$
as 
\begin{equation}
D(\sigma_{y})=|\{(i,j):y_{j}<y_{i}\mbox{ and }i<j\}|.\label{alternativeD}
\end{equation}
Given a vector $y$ and indices $k$ and $r$, let
\[
D_{k,r}(y)=|\{(i,j):y_{j}<y_{i}\mbox{ and }i+k\leq j\mbox{ and } r\leq i\}|.
\]
In particular, $D_{1,1}(y)=D(\sigma_{y})$.
With this definition, the metric $D_{k,r}$ enables us to quantify inverted pairs locally inside specific ranges of the set $y$, while Kendall distance only provides a global measure. This way, it becomes more interesting for our purposes, as the next Theorem illustrates.
\begin{thm}
\label{thm:maxJump}Let $y$ be the second eigenvector of the random matrix $\widehat{M}$. Let $r=n^{\alpha}$ and $k=n^{\beta}$. Then with probability $1-n^{-3}$ there is a constant $C>0$, such that 
if $\alpha\leq \beta$, then $D_{k,r}(y)< Cn^{5-2(\alpha+\beta)}$, and
if $\alpha > \beta$, then $D_{k,r}(y)< Cn^{5-4\beta}$.
\end{thm}
Here, whenever $\alpha+\beta > 3/2$ and $\alpha\leq \beta$, or whenever $\beta>3/4$ and $\alpha < \beta$, we get the extremal cases where the number of inverted pairs satisfy $D_{k,r}(y)\in o(n^2)$.

\begin{thm}
\label{thm:maxJump2}Let $y$ be the second eigenvector of the random matrix $\widehat{M}$, and let $k=n^{\beta}$. Then with probability $1-n^{-3}$ there is a constant $C>0$, such that 
 $D_{k,1}(y)< Cn^{\frac{7-2 \beta}{3}}$.
\end{thm}

Now, if we take $\beta=1/2$, then $D_{k,1}\in o(n^2)$, which means almost no vertex will drift more than $\sqrt{n}$ from its correct position. That ressembles the degree approach explained above, where the position of
most vertices will drift at most $O(\sqrt{n} \log n)$ positions from their correct position.

Theorems \ref{thm:maxJump} and \ref{thm:maxJump2} expose an interesting behavior in the group of bad
vertices. There is a trade off between how far vertices can jump and
the total number of such incorrectly placed vertices. That is useful
for our purpose to establish metrics on the correctness of the rank. The next Theorem shows that the permutation derived from the second
eigenvector $\hat{x}$ of the random matrix $\widehat{M}$ is well
behaved.
\begin{thm}
\label{thm:footrule}Let $y$ be the second eigenvector of the random
matrix $\widehat{M}$.
Then $D(\sigma_y)=\mathcal{O}(n^{9/5})$ with probability $1-n^{-3}$. 

\begin{comment}
\begin{thm}
Let $\sigma$ be the permutation derived from the second eigenvector
$\hat{x}$ of the random matrix $\widehat{M}$, then $F(\sigma)\leq\mathcal{O}(n^{7/4})$
with probability $1-n^{-3}$. \end{thm}
\end{comment}

\end{thm}
As a consequence, Kendall $\tau$ rank correlation coefficient is 
\begin{cor}
$\tau=1-\mathcal{O}(n^{-1/5})$ with probability $1-n^{-3}$.
\end{cor}
To prove Theorem \ref{thm:upperBound} and quantify the error in the
order given by Algorithm \ref{alg:algo1}, our technique strongly
relies on analytic expressions for the eigenvalues and eigenvectors
of the model matrix. To this end, the next section provides the necessary
information about the spectrum of the model matrix. One of the outcomes
of Section \ref{sec:eigenvalues} is that one eigenvector of the model
matrix is steep enough. That means it contains the desired information
about the structure of the graph, and provides the correct order of
vertices. Moreover, consecutive entries differ significantly enough so that a small perturbation will have limited effect on the order. Further, in Section \ref{sec:eigRandom} we show that this
steep eigenvector is close in norm to the eigenvector of the random
graph. Thus, we prove Theorem \ref{thm:upperBound} in Section \ref{sec:eigRandom}.
In Section \ref{sec:qualitative} we perform a qualitative analysis
of the problem and prove Theorems \ref{thm:maxJump} and \ref{thm:footrule}.

\section{\label{sec:eigenvalues}The eigenvalues of the model matrix}

In this section we look to the spectrum of the model matrix which
we identify by $M$. Consider the band matrix of even order $n$ defined
as{\tiny{}
\[
A:=\left[\begin{array}{cccccccc}
1 & 1 & \ldots & 1 & 0 & 0 & \ldots & 0\\
1 & 1 &  &  & 1 &  &  & 0\\
\vdots & \vdots & \ddots &  &  & \ddots &  & \vdots\\
1 & 1 &  &  &  &  & 1 & 0\\
0 & 1 &  &  &  &  &  & 1\\
0 &  & \ddots &  &  &  &  & 1\\
\vdots &  &  &  &  &  &  & \vdots\\
0 & 0 & \ldots & 0 & 1 & 1 & \ldots & 1
\end{array}\right],
\]
}where exactly $s:=n/2$ columns have the first entry equal to $1$.
Notice that $M=p(A-I)$. Thus $M$ and $A$ share the same eigenvectors
with eigenvalues related as $\lambda(M)=p(\lambda(A)-1)$. For simplicity
we proceed investigating the spectrum of $A$. 

This matrix is an example of the familiar Toeplitz matrix which arises in
many applications. Though Toeplitz matrices are well studied, their
spectrum is unknown in the general form \cite{boettcher,Gray}. Finding the spectrum
seems to be a hard problem. There have been some advances in finding
analytic expressions for the spectrum of particular instances, but
to the best of our knowledge the spectrum of the matrix $A$ is unknown.
In this section we find expressions for the $s$ eigenvalues and eigenvectors
of $A$. Among them we identify the second largest, in absolute value,
which is the main tool in this paper. We also characterize the remaining
eigenvalues of $A$, and approximate its largest eigenvalue. We will need this information to bound the error between eigenvectors of model matrix and random matrix. 

To this end, label the eigenvalues of $A$ as 
\[
\left|\lambda_{1}\left(A\right)\right|\geq\left|\lambda_{2}\left(A\right)\right|\geq\ldots\geq\left|\lambda_{n}\left(A\right)\right|.
\]

In this section we prove the following result, which shows that the
entries of the second eigenvector of the model matrix are monotonic.
\begin{thm}
\label{thm:eigA}The eigenvalue $\left|\lambda_{2}\left(A\right)\right|=1/\sqrt{2+2\cos\left(\frac{2s\pi}{2s+1}\right)}$
with corresponding eigenvector $u=\left[u_{j}\right]$ where $u_{j}=\cos\left(\frac{(2j-1)\pi}{4s+2}\right)$,
for $j=1,\ldots,s$ and $u_{j}=-u_{n-j+1}$, for $j=s+1,\ldots,n$,
where $s=n/2$.
\end{thm}
In order to prove Theorem \ref{thm:eigA} we need to compute the eigenvalues
of auxiliary matrices. Let $J$ be the matrix with all entries equal
to $1$ and let $B:=J-A$. Then $B$ can be written in the form {\tiny{}
\[
B:=\left[\begin{array}{cccccccc}
0 & 0 & \ldots & 0 & 1 & 1 & \ldots & 1\\
0 & 0 &  &  & 0 &  &  & 1\\
\vdots & \vdots & \ddots &  &  & \ddots &  & \vdots\\
0 & 0 &  &  &  &  & 0 & 1\\
1 & 0 &  &  &  &  &  & 0\\
1 &  & \ddots &  &  &  &  & 0\\
\vdots &  &  &  &  &  &  & \vdots\\
1 & 1 & \ldots & 1 & 0 & 0 & \ldots & 0
\end{array}\right]=\left[\begin{array}{cc}
 & C^{T}\\
C
\end{array}\right],\text{ where }C:=\left[\begin{array}{cccc}
1 & 0 & \cdots & 0\\
1 & 1 &  & \vdots\\
\vdots &  & \ddots & 0\\
1 & \ldots & 1 & 1
\end{array}\right]
\]
}is the lower triangular matrix of order $s=n/2$. Define the permutation
matrix {\tiny{}
\[
P:=\left[\begin{array}{cccc}
 &  &  & 1\\
 &  & 1\\
 & 1\\
1
\end{array}\right].
\]
}{\tiny \par}
\begin{lem}
\label{lem:CP}Let $x$ be any eigenvector for the matrix $CP$ with
eigenvalue $\lambda.$ Define vectors $y:=\left[\begin{array}{cc}
Px^{T} & -x^{T}\end{array}\right]^{T}$ and $z:=\left[\begin{array}{cc}
Px^{T} & x^{T}\end{array}\right]^{T}$. Then $By=-\lambda y$ and $Bz=\lambda z$.\end{lem}
\begin{proof}
Now for any eigenvector $x$ for the matrix $CP$ we have $CPx=\lambda x$
and $C^{T}x=\lambda Px$. Thus we obtain
\[
By=\left[\begin{array}{c}
-C^{T}x\\
CPx
\end{array}\right]=\left[\begin{array}{c}
-\lambda Px\\
\lambda x
\end{array}\right]=-\lambda y.
\]
In the same way we can see that $Bz=\lambda z$. This finishes the
proof.
\end{proof}
In light of Lemma \ref{lem:CP}, it suffice to compute the eigenvalues
of $CP$ in order to get the eigenvalues of $B$. The next lemma describes
these eigenvalues and its eigenvectors.
\begin{thm}
\label{thm:eigB}For $k=1\ldots s$, the eigenvalues of $B$ are $\lambda_{k}=1/\sqrt{2+2\cos\left(\frac{2k\pi}{2s+1}\right)}$
with corresponding eigenvectors $u^{k}=\left[u_{j}\right]$ where
$u_{j}=\left(-1\right)^{j}\sin\left(\frac{2jk\pi}{2s+1}\right)$,
for $j=1,\ldots,s$ and $u_{j}=-u_{n-j+1}$, for $j=s+1,\ldots,n$. 

For $k=s+1\ldots n$ we have $\lambda_{k}=-1/\sqrt{2+2\cos\left(\frac{2k\pi}{2s+1}\right)}$
with corresponding eigenvectors $u^{k}=\left[u_{j}\right]$ where
$u_{j}=\left(-1\right)^{j}\sin\left(\frac{2jk\pi}{2s+1}\right)$,
for $j=1,\ldots,s$ and $u_{j}=u_{n-j+1}$, for $j=s+1,\ldots,n$. \end{thm}
\begin{proof}
We will show that for $k=1\ldots s$ the eigenvalues $\lambda_{k}$
of $CP$ are $1/\sqrt{2+2\cos\left(\frac{2k\pi}{2s+1}\right)}$ with
corresponding eigenvectors $u^{k}=\left[u_{j}\right]$ where $u_{j}=\left(-1\right)^{j}\sin\left(\frac{2jk\pi}{2s+1}\right)$,
for $j=1,\ldots,s$. Then, by Lemma \ref{lem:CP} the proof is done.

We notice that {\tiny{}
\[
CP=\left[\begin{array}{cccc}
0 & \cdots & 0 & 1\\
0 &  & \iddots & \vdots\\
\vdots & 1 &  & 1\\
1 & \ldots & 1 & 1
\end{array}\right],\left(CP\right)^{2}=\left[\begin{array}{ccccc}
1 & 1 & \cdots & 1 & 1\\
1 & 2 &  & 2 & 2\\
\vdots &  & \ddots & \vdots & \vdots\\
1 & 2 & \cdots & s-1 & s-1\\
1 & 2 & \cdots & s-1 & s
\end{array}\right]\text{, and }\left(CP\right)^{-2}=\left[\begin{array}{ccccc}
2 & -1\\
-1 & 2 & -1\\
 & -1 & \ddots & \ddots\\
 &  & \ddots & 2 & -1\\
 &  &  & -1 & 1
\end{array}\right].
\]
}{\tiny \par}

That is familiar tridiagonal matrix which appears often in the literature.
One of the oldest works where its eigenvalues were reported is in
\cite{severalTridiag}. We finish the proof by indicating that the
tridiagonal matrix $\left(CP\right)^{-2}$ has eigenvalues $2+2\cos\left(\frac{2k\pi}{2s+1}\right)$
for $k=1\ldots s$, with corresponding eigenvector $u=\left[u_{j}\right]$
where $u_{j}=\left(-1\right)^{j}\sin\left(\frac{2jk\pi}{2s+1}\right)$,
for $j=1,\ldots,s$. Thus, $u$ is an eigenvector for $CP$ with corresponding
eigenvalue $\lambda_{k}=1/\sqrt{2+2\cos\left(\frac{2k\pi}{2s+1}\right)}$.
This concludes the proof.
\end{proof}
Now, we are ready to prove the following theorem concerning the eigenvalues
of $A$.
\begin{thm}
\label{thm:AandBsameEig}$\left|\lambda_{2}\left(A\right)\right|=\left|\gamma_{2}\right|$,
where $\gamma_{2}$ is the second largest eigenvalue of $B$ in absolute
value. Besides, both eigenvalues have the same eigenvector.\end{thm}
\begin{proof}
First, we notice that all matrices $A,$ $B$ and $J$ satisfy the Perron-Frobenius
Theorem. Therefore, they all have a unique eigenvector with entries
of the same sign corresponding to the largest eigenvalue in absolute
value. They are called Perron vector and Perron value for the matrix
of interest. Next, we define the set $H$ which is just the set of
unit vectors excluding all possible Perron vectors. 
\[
H:=\left\{ x\in\mathbb{R}^{n}\backslash\left\{ \mathbb{R}_{+}^{n}\cup\mathbb{R}_{-}^{n}\right\} :\left\Vert x\right\Vert =1\right\} .
\]

The Perron value of $A$ can be expressed as 
\[
\left|\lambda_{1}(A)\right|=\max_{\left\Vert x\right\Vert =1}x^{T}Ax.
\]
Thus, if we restrict the set of vectors where we take the maximum
to the set $H$, we get 
\[
\left|\lambda_{2}(A)\right|=\max_{x\in H}\left|x^{T}Ax\right|.
\]
And the same maximization problem can be used to identify $\left|\gamma_{2}\right|$
and $\left|\lambda_{2}(J)\right|$. Finally, we can write

\begin{eqnarray*}
\left|\lambda_{2}(A)\right| & = & \max_{x\in H}\left|x^{T}(J-B)x\right|\\
 & \leq & \max_{x\in H}\left\{ \left|x^{T}Jx\right|+\left|x^{T}Bx\right|\right\} \\
 & \leq & \max_{x\in H}\left|x^{T}Jx\right|+\max_{x\in H}\left|x^{T}Bx\right|\\
 & = & \left|\lambda_{2}(J)\right|+\left|\gamma_{2}\right|.
\end{eqnarray*}
Besides, the only nonzero eigenvalue of $J$ is its Perron value.
Therefore, the above inequality is simply 
\[
\left|\lambda_{2}(A)\right|\leq\left|\gamma_{2}\right|.
\]
By noticing that $B=J-A$, the same calculation gives us the inequality
$\left|\gamma_{2}\right|\leq\left|\lambda_{2}(A)\right|.$ Both inequalities
prove the first statement of the theorem. 

To see that $\lambda_{2}(A)$ and $\gamma_{2}$ have the same eigenvector
we apply Theorem \ref{thm:eigB}. According to this theorem, there is an eigenvector $u$
of $\gamma_{2}$, with positive and negative entries, such that $u=\left[\begin{array}{cc}
Pv^{T} & -v^{T}\end{array}\right]^{T}$. Clearly $Ju=0$ and then 
\[
Au=(J-B)u=-Bu=-\gamma_{2}u.
\]
That finishes the proof.
\end{proof}
Next, we prove the main theorem of this section.
\begin{proof}
(Theorem \ref{thm:eigA}) Now the proof of Theorem \ref{thm:eigA}
follows easily from Theorem \ref{thm:eigB} and Theorem \ref{thm:AandBsameEig}.
The largest eigenvalue of $B$ occurs for $k=1$ in Theorem \ref{thm:eigB}.
Relabelling the eigenvalues $\lambda_{i}(B)$ for $i=1\ldots n$ such
that$\left|\gamma_{1}\right|\geq\left|\gamma_{2}\right|\geq\ldots\geq\left|\gamma_{n}\right|$,
we see that $\left|\gamma_{1}\right|=\left|\gamma_{2}\right|$. Consider
$k=1$ in Theorem \ref{thm:eigB}. The entries of the eigenvector
of $\lambda_{1}(A)$ can be rewritten as $u_{j}=\cos\left(\frac{(2j-1)\pi}{4s+2}\right)$,
for $j=1,\ldots,s$, and that establishes Theorem \ref{thm:eigA}.
\end{proof}

\section{\label{sec:eigRandom}The eigenvectors of the random graph}

We apply random matrix theory and the results about the spectrum of
the model matrix, obtained above, to show that the second eigenvector
of the random graph is, with high probability, close to the second
eigenvector of the model. Thus the second eigenvector of the random
matrix can be used to reveal the correct linear structure of the random
graph. 

In the last few years, there has been significant progress on the
subject of random matrices, especially around the universality phenomenon
\cite{TaoVu-Universality-1} and concentration inequalities (see \cite{Tao,Vu}
for good surveys). For our purposes, we will use the following concentration
inequality from \cite{Vu-1}.
\begin{lem}
\label{lem:normRand}(Norm of a random matrix). There is a constant
$C>0$ such that the following holds. Let $E$ be a symmetric matrix
whose upper diagonal entries $e_{ij}$ are independent random variables
where $e_{ij}=1-p_{ij}$ or $-p_{ij}$ with probabilities $p_{ij}$
and $1-p_{ij}$ , respectively, where $0\leq p_{ij}\leq1$. Let $\sigma^{2}=\max_{ij}p_{ij}(1-p_{ij})$.
If $\sigma^{2}\geq C\log n/n$, then 
\[
\mathbb{P}(\left\Vert E\right\Vert \geq C\sigma n^{1/2})\leq n^{-3}.
\]

\end{lem}
A well-known result by Davis and Kahan \cite{Davis-Kahan}, from matrix
theory says that the angle between eigenvectors of the model matrix
and of the random one is bounded in terms of their spectrum:
\begin{lem}
\label{lem:Davis-Kahan} (Davis-Kahan) Let $A$ and $\widehat{A}$
be symmetric matrices, $\lambda_{i}$ be the $i-$th largest eigenvalue
of $A$ with eigenvector $x_{i}$, and $\widehat{x_{i}}$ be the eigenvector
of the $i-$th largest eigenvalue of $\widehat{A}$. If $\theta$
is the angle between $x_{i}$ and $\widehat{x_{i}}$, then 
\[
\sin2\theta\leq\frac{2\left\Vert A-\widehat{A}\right\Vert }{\min_{i\neq j}\left|\lambda_{i}-\lambda_{j}\right|},\text{ provided }\lambda_{i}\neq\lambda_{j}.
\]

\end{lem}
To apply the Davis-Kahan Theorem to $\lambda_{2}$, we need the smallest
gap between the eigenvalues. The next Theorem helps us with that.
\begin{thm}
\label{thm:normA}$\lambda_{1}(A)=1/\sqrt{2-2\cos\left(\theta_{1}\right)}$,
where $\theta_{1}\leq\frac{\pi}{4s}$ and $\lambda_{3}(A)=1/\sqrt{2-2\cos\left(\theta_{3}\right)}$,
where $\theta_{3}>\frac{\pi}{s}$.
\end{thm}
To prove Theorem \ref{thm:normA} and characterize the additional
eigenvalues of $A$ we consider the following matrix 

{\tiny{}
\[
D=\left[\begin{array}{ccccc}
1 & \cdots &  & 1 & 1\\
1 &  &  & \iddots & 2\\
 &  &  & \iddots & \vdots\\
\vdots & 1 &  & 2 & 2\\
1 & 2 & \ldots & 2 & 2
\end{array}\right].
\]
}{\tiny \par}
\begin{lem}
\label{lem:eigDandA}Let $x$ be any eigenvector for matrix $D$ with
eigenvalue $\lambda.$ Define the vector $y:=\left[\begin{array}{cc}
x^{T} & Px^{T}\end{array}\right]^{T}$. Then $Ay=\lambda y$.\end{lem}
\begin{proof}
First, we have{\tiny{} $A\left[\begin{array}{c}
I\\
P
\end{array}\right]=\left[\begin{array}{cc}
J & C\\
C^{T} & J
\end{array}\right]\left[\begin{array}{c}
I\\
P
\end{array}\right]=\left[\begin{array}{c}
J+CP\\
C^{T}+JP
\end{array}\right]=\left[\begin{array}{c}
D\\
PD
\end{array}\right].$} Then{\tiny{} $Ay=A\left[\begin{array}{c}
I\\
P
\end{array}\right]x=\left[\begin{array}{c}
D\\
PD
\end{array}\right]x=\lambda\left[\begin{array}{c}
x\\
Px
\end{array}\right]=\lambda y.$}{\tiny \par}
\end{proof}
The previous lemma leads us to the study of the spectrum of the matrix $D$.
\begin{thm}
\label{thm:eigenvaluesD}Let $D$ be as defined above, of even order $s$.
Then, its eigenvalues are $\lambda_{k}=1/\sqrt{2-2\cos\left(\theta_k\right)}$,
for $k=1\ldots s$, where $\theta_k$ is a root of the equation 
\[
p(\theta)=\sin\left(\left(s+1\right)\theta\right)+3\sin\left(s\theta\right)-4\sin\left(\left(s-1\right)\theta\right)-4\sin\left(\theta\right).
\]
\end{thm}
\begin{proof}
It is easy to verify that {\tiny{}
\[
D^{-1}=\left[\begin{array}{ccccc}
2 & 0 & \cdots &  & -1\\
0 &  &  & -1 & 1\\
\vdots &  & \iddots & \iddots & 0\\
 &  &  &  & \vdots\\
-1 & 1 & 0 & \cdots & 0
\end{array}\right]\text{ and }D^{-2}=\left[\begin{array}{ccccc}
5 & -1 &  &  & -2\\
-1 & 2 & -1\\
 & -1 & \ddots & \ddots\\
 &  & \ddots & 2 & -1\\
-2 &  &  & -1 & 2
\end{array}\right].
\]
} Therefore, it suffices to prove that the eigenvalues of $D^{-2}$
are of the required form $2-2\cos\left(\theta\right)$.

In order to do that, we obtain an expression for the characteristic
polynomial of $D^{-2}$. A relevant matrix related to this expression
is the order $s$ matrix defined as {\tiny{}
\[
E_{s}=\left[\begin{array}{ccccc}
2-\lambda & -1 & 0 & \cdots & 0\\
-1 & 2-\lambda &  & \ddots & \vdots\\
0 &  &  &  & 0\\
\vdots & \ddots &  & \ddots & -1\\
0 & \cdots & 0 & -1 & 2-\lambda
\end{array}\right].
\]
}Below $\left|M\right|_{n}$ indicates the determinant of the order
$n$ matrix $M$. By expanding the first row the determinant of $D^{-2}-\lambda I$,
we can write{\tiny{}
\begin{eqnarray*}
det(D^{-2}-\lambda I) & = & (5-\lambda)\left|E_{s-1}\right|\\
 &  & -\left|\begin{array}{ccccc}
-1 & -1 & 0 & \cdots & 0\\
0 & 2-\lambda & -1 & 0\\
\vdots & -1 & \ddots & \ddots\\
0 &  & \ddots &  & -1\\
-2 & 0 & \cdots & -1 & 2-\lambda
\end{array}\right|_{s-1}+2\left|\begin{array}{ccccc}
-1 & 2-\lambda & -1 &  & 0\\
0 & -1 & 2-\lambda & \ddots\\
\vdots & 0 & \ddots & \ddots & -1\\
0 &  & \ddots & -1 & 2-\lambda\\
-2 & 0 & \cdots & 0 & -1
\end{array}\right|_{s-1}.
\end{eqnarray*}
}Then we expand the second and third determinants by the first column,
which gives us

{\tiny{}
\begin{eqnarray*}
det(D^{-2}-\lambda I) & = & (5-\lambda)\left|E_{s-1}\right|\\
 &  & -\left(\left|E_{s-2}\right|+2\left|\begin{array}{ccccc}
-1 & 0 & \cdots &  & 0\\
2-\lambda & -1 & 0 &  & 0\\
-1 & \ddots & \ddots & \ddots & \vdots\\
 & \ddots &  &  & 0\\
0 &  & -1 & 2-\lambda & -1
\end{array}\right|_{s-2}\right)\\
 &  & +2\left(-\left|\begin{array}{ccccc}
-1 & 2-\lambda & -1 &  & 0\\
0 & -1 & 2-\lambda & \ddots\\
 & 0 & \ddots & \ddots\\
 &  & \ddots\\
 &  & \ddots & -1 & 2-\lambda\\
0 &  & \cdots & 0 & -1
\end{array}\right|_{s-2}-2\left|E_{s-2}\right|\right)\\
 & = & (5-x)\left|E_{s-1}\right|-5\left|E_{s-2}\right|.\\
 &  & -4\left|\begin{array}{ccccc}
-1 & 2-\lambda & -1 &  & 0\\
0 & -1 & 2-\lambda & \ddots\\
 & 0 & \ddots & \ddots\\
 &  & \ddots\\
 &  & \ddots & -1 & 2-\lambda\\
0 &  & \cdots & 0 & -1
\end{array}\right|_{s-2}
\end{eqnarray*}
}Since the last matrix has even order and $-1$ in the diagonal, we
get

\begin{equation}
det(D^{-2}-\lambda I)=(5-\lambda)\left|E_{s-1}\right|-5\left|E_{s-2}\right|-4.\label{eq:detD-2}
\end{equation}

From this point, we call up Chebyshev polynomials of second kind.
They are defined by the recurrence relation \begin{align*} U_0(x) & = 1 \\ U_1(x) & = 2x \\ U_{n+1}(x) & = 2xU_n(x) - U_{n-1}(x). \end{align*} Besides,
it is a well-known fact that $U_{n}(x)$ obeys the determinant identity{\tiny{}
\[
U_{n}(x)=\left|\begin{array}{cccc}
2x & -1 & 0 & 0\\
-1 & 2x & \ddots & 0\\
0 & \ddots & \ddots & -1\\
0 & 0 & -1 & 2x
\end{array}\right|_{n}.
\]
}Therefore, a change of variable $x=\frac{2-\lambda}{2}$ provides
us the determinant 
\[
\left|E_{s}\right|=U_{s}\left(\frac{2-\lambda}{2}\right)
\]
and the relation 
\[
U_{s+1}\left(\frac{2-\lambda}{2}\right)=\left(2-\lambda\right)U_{s}\left(\frac{2-\lambda}{2}\right)-U_{s-1}\left(\frac{2-\lambda}{2}\right),
\]
which we apply in the equivalent form
\[
\left(2-\lambda\right)U_{s}\left(\frac{2-\lambda}{2}\right)=U_{s+1}\left(\frac{2-\lambda}{2}\right)+U_{s-1}\left(\frac{2-\lambda}{2}\right).
\]
In this manner, equation (\ref{eq:detD-2}) becomes

\begin{eqnarray*}
det(D^{-2}-\lambda I) & = & (3+2-\lambda)U_{s-1}\left(\frac{2-\lambda}{2}\right)-5U_{s-2}\left(\frac{2-\lambda}{2}\right)-4\\
 & = & U_{s}\left(\frac{2-\lambda}{2}\right)+3U_{s-1}\left(\frac{2-\lambda}{2}\right)-4U_{s-2}\left(\frac{2-\lambda}{2}\right)-4.
\end{eqnarray*}

Furthermore, one of many well-known properties of Chebyshev polynomials
is 
\[
U_{n}(cos\theta)=\frac{\sin\left(\left(n+1\right)\theta\right)}{\sin\theta}.
\]
Therefore, the change of variables $\cos\theta=\frac{2-\lambda}{2}$
allows us to rewrite $det(D^{-2}-\lambda I)$ as 
\[
det(D^{-2}-(2-2\cos\theta)I)=\frac{\sin\left(\left(s+1\right)\theta\right)+3\sin\left(s\theta\right)-4\sin\left(\left(s-1\right)\theta\right)-4\sin\left(\theta\right)}{\sin\left(\theta\right)}.
\]
Finally, $2-2\cos\theta$ is an eigenvalue of $D^{-2}$ whenever $\theta$
is a root of the equation 
\[
p(\theta)=\sin\left(\left(s+1\right)\theta\right)+3\sin\left(s\theta\right)-4\sin\left(\left(s-1\right)\theta\right)-4\sin\left(\theta\right),
\]
which concludes the proof.\end{proof}
\begin{thm}
\label{thm:boundEigD}Let $D$ be a matrix as defined earlier, of even order $s\geq 5$.
Then, its largest eigenvalues 
\[
\lambda_{i}=1/\sqrt{2-2\cos\left(\theta_{i}\right)},
\]
for $i=1,2$ are such that $\theta_{1}\in\left(\frac{\pi}{s^{2}},\frac{\pi}{4s}\right)$
and $\theta_{2}>\frac{\pi}{s}$. \end{thm}
\begin{proof}
To bound $\theta_{1}$, Theorem \ref{thm:eigenvaluesD} guarantees
it is enough to locate the smallest root of 
\[
p(\theta)=\sin\left(\left(s+1\right)\theta\right)+3\sin\left(s\theta\right)-4\sin\left(\left(s-1\right)\theta\right)-4\sin\left(\theta\right).
\]
On one hand, by means of trigonometric identities we get 
\[
p\left(\frac{x\pi}{s}\right)=-3\,\sin\left(x\pi\right)\cos\left({\frac{x\pi}{s}}\right)+5\,\cos\left(x\pi\right)\sin\left({\frac{x\pi}{s}}\right)+3\,\sin\left(x\pi\right)-4\,\sin\left({\frac{x\pi}{s}}\right).
\]
We claim that $p\left(\frac{x\pi}{s}\right)<0$ for $x\in[1/4,1)$.
To see it, it is enough to show that 
\[
{\frac{5\,\cos\left(x\pi\right)-4}{3\sin\left(x\pi\right)}}<\frac{\cos\left(\frac{x\pi}{s}\right)-1}{\sin\left(\frac{x\pi}{s}\right)},
\]
for $x\in[1/4,1)$. Notice that $\frac{5\,\cos\left(x\pi\right)-4}{3\sin\left(x\pi\right)}$
is a decreasing function in $x$ for $x\in[1/4,1)$. Thus, it is
enough to show that
\[
{\frac{5\,\cos\left(\pi/4\right)-4}{3\sin\left(\pi/4\right)}}<\frac{\cos\left(\frac{x\pi}{s}\right)-1}{\sin\left(\frac{x\pi}{s}\right)}.
\]
Also $\frac{\cos\left(x\frac{\pi}{s}\right)-1}{\sin\left(x\frac{\pi}{s}\right)}$
is an decreasing function in $x$ for $x\in[1/4,1]$ and increasing
as a function in $s$. Thus, it is enough to check that
\[
{\frac{5\,\cos\left(\pi/4\right)-4}{3\sin\left(\pi/4\right)}}<\frac{\cos\left(\frac{\pi}{s}\right)-1}{\sin\left(\frac{\pi}{s}\right)}.
\]
In fact, the inequality is true for $s\geq 2$, thus the claim follows.

On the other hand, the Taylor series of $p\left(\theta\right)$ at
$\theta=0$ is

\[
\theta+\left(-\frac{5}{2}\,s^{2}+\frac{3}{2}\,s-\frac{1}{6}\right)\theta^{3}+\left({\frac{5}{24}}\,{s}^{4}-\frac{1}{4}\,{s}^{3}+{\frac{5}{12}}\,{s}^{2}-\frac{1}{8}\,s+{\frac{1}{120}}\right){\theta}^{5}+\mathcal{O}\left({\theta}^{6}\right).
\]
Furthermore, if $s$ is not too small the terms of order 5 or greater
adds up to a positive constant, thus 
\[
f(\theta)=\theta+\left(-\frac{5}{2}\,s^{2}+\frac{3}{2}\,s-\frac{1}{6}\right)\theta^{3}\leq p\left(\theta\right)
\]
and 
\[
f\left(\frac{\pi}{s^{2}}\right)=\frac{\pi\left(6s^{4}-15\pi^{2}s^{2}+9\pi^{2}s-\pi^{2}\right)}{6s^{6}}.
\]
Now, the polynomial $6s^{4}-15\pi^{2}s^{2}+9\pi^{2}s-\pi^{2}$ in
the variable $s$ has largest root smaller than $5$. Therefore,
\[
0\leq f\left(\frac{\pi}{s^{2}}\right)\leq p\left(\frac{\pi}{s^{2}}\right)\text{ and }p\left(\frac{\pi}{4s}\right)<0.
\]
Then, by continuity $p\left(\theta\right)$ has a root in $\left(\frac{\pi}{s^{2}},\frac{\pi}{4s}\right)$.
Finally, the fact that $p\left(\frac{x\pi}{s}\right)<0$ for $x\in[1/4,1)$
guarantees that $\theta_{2}>\frac{\pi}{s}$, which concludes the proof.
\end{proof}
Finally, we are able to look at $\lambda_{1}(A)$ and $\lambda_{3}(A)$
and prove Theorem \ref{thm:normA}.
\begin{proof}
(Theorem \ref{thm:normA}) By Lemma \ref{lem:eigDandA}, each eigenvalue
of $D$ is an eigenvalue of $A$. Thus, by Theorem \ref{thm:boundEigD}
the proof is done for $\lambda_{1}$ if we prove that $\lambda_{1}(A)=\lambda_{1}(D)$.
But that is easy to see, since $D$ and $A$ are matrices that fulfill
Perron-Frobenius Theorem. Thus, its unique largest eigenvalue has
an eigenvector with entries of the same sign, so called Perron vector.
Again by Lemma \ref{lem:eigDandA}, if $x$ is a Perron vector of
$D$, then $\left[\begin{array}{cc}
x^{T} & Px^{T}\end{array}\right]^{T}$ is an eigenvector of $A$ with entries of the same sign. Thus it
must be a Perron vector of $A$, which gives us $\lambda_{1}(A)=\lambda_{1}(D)$. 

Now $\lambda_{3}(A)$ is characterized by Lemma \ref{lem:eigDandA}, which provides the eigenvalues with odd indices for $A$,
and Theorem \ref{thm:boundEigD}. Thus $\lambda_{3}(A)=1/\sqrt{2-2\cos\left(\theta_{3}\right)}$,
where $\theta_{3}>\frac{\pi}{s}$. That finishes the proof.
\end{proof}
Finally, the results of this section enables us to prove Theorem \ref{thm:upperBound}.
\begin{proof}
(Theorem \ref{thm:upperBound}) We can view the adjacency matrix $\widehat{M}$
as a perturbation of $M$, $\widehat{M}=M+E$, where the entries of
$E$ are $e_{ij}=1-p$ with probability $p$ and $-p$ with probability
$1-p$. That way $E$ is as in Lemma \ref{lem:normRand}, and then
with probability at least $1-n^{-3}$, we have 

\begin{equation}
\left\Vert E\right\Vert \leq C\sigma\sqrt{n},\label{eq:difference}
\end{equation}
for some constant $C>0$. Since $x$ and $\widehat{x}$ are both unitary
we have $\left\Vert x-\widehat{x}\right\Vert \leq\sqrt{2}\sin\theta$,
where $\theta$ is the angle they form. Therefore, we can apply Lemma
\ref{lem:Davis-Kahan} and equation (\ref{eq:difference}), to obtain
a constant $C_{1}>0$ such that

\begin{eqnarray}
\left\Vert x-\widehat{x}\right\Vert  & \leq & \frac{C_{1}\left\Vert M-\widehat{M}\right\Vert }{\min_{i\neq2}\left|\lambda_{i}-\lambda_{2}\right|},\nonumber \\
 & \leq & \frac{C_{1}\sigma\sqrt{n}}{\min_{i\neq2}\left|\lambda_{i}-\lambda_{2}\right|}.\label{eq:boundsss}
\end{eqnarray}

Also, if $A$ is as in Theorem \ref{thm:eigA}, we have $M=p(A-I).$
Now, to bound the gap between the eigenvalues, we apply Theorems \ref{thm:eigA}
and \ref{thm:normA}, to get

\[
\lambda_{1}(M)=p\lambda_{1}(A-I)=p/\sqrt{2-2\cos\left(t_{1}\right)}-p,
\]
\[
\lambda_{3}(M)=p\lambda_{3}(A-I)=p/\sqrt{2-2\cos\left(t_{3}\right)}-p\text{ and}
\]
\[
\lambda_{2}(M)=p\lambda_{2}(A-I)=p/\sqrt{2+2\cos\left(\frac{2s\pi}{2s+1}\right)}-p,
\]
where $t_{1}\leq\frac{\pi}{4s}$, $t_{3}>\frac{\pi}{s}$, and $s=\frac{n}{2}$
. Notice that, since $\cos\left(\theta_{1}\right)$ is decreasing
in $\theta_{1}$, we have 
\[
\frac{1}{\lambda_{1}-\lambda_{2}}\leq\frac{1}{p/\sqrt{2-2\cos\left(\frac{\pi}{4s}\right)}-p/\sqrt{2+2\cos\left(\frac{2s\pi}{2s+1}\right)}}.
\]
Now, an asymptotic expansion of the last expression, shows that there
is a constant $B_1>0$ such that 
\[
\frac{1}{\lambda_{1}-\lambda_{2}}\leq\frac{B_1}{s}.
\]
Similarly, we obtain a constant $B_2>0$ such that 
\[
\frac{1}{\lambda_{3}-\lambda_{1}}\leq\frac{B_2}{s}.
\]

Therefore, for some constant $C_{2}>0$ inequality \ref{eq:boundsss}
becomes
\[
\left\Vert x-\widehat{x}\right\Vert \leq\frac{C_{2}\sqrt{n}}{s}=\frac{2C_{2}\sqrt{n}}{n}.
\]
That finishes the proof.
\end{proof}

\section{\label{sec:qualitative}Bounding the number of misplaced vertices}

Our method relies on the asymptotic expansions of the terms $x_{i}$
regarding it as a function of $n$, as the next Lemma states.
\begin{lem}
\label{lem:Xunitary}Let $x$ be the unitary eigenvector for $\lambda_{2}\left(A\right)$.
Then 
\[
\left(x_{r}-x_{r+k}\right)^{2}=\frac{1}{\pi}\,k^{2}\left(2\,r+k-1\right)^{2}\theta^{5}+O\left(\theta^{7}\right),
\]
where $\theta=\frac{\pi}{2s+1}$.\end{lem}
\begin{proof}
Theorem \ref{thm:eigA} provides the expressions we need to compute
$x_{j}$: we have $x_{j}=\omega\cos\left(\frac{(2j-1)\pi}{4s+2}\right)$,
for $j=1,\ldots,s$, where $\omega$ is a constant such that $\left\Vert x\right\Vert =1$.
Thus, 
\[
\left(x_{r}-x_{r+k}\right)^{2}=\omega^{2}\left(\cos\left(\frac{(2r-1)\pi}{4s+2}\right)-\cos\left(\frac{(2(r+k)-1)\pi}{4s+2}\right)\right)^{2}.
\]
To find $\omega^{2}$, we make use of the trigonometric identity 
\[
\sum_{i=0}^{k-1}\cos(2\alpha i+\beta)=\frac{\sin\left(k\alpha\right)\cos\left(\beta+\left(k-1\right)\alpha\right)}{\sin\alpha}.
\]
Then, we can write

\begin{eqnarray}
c & = & \sum_{i=1}^{s}\cos^{2}\left(\frac{(2i-1)\pi}{4s+2}\right)\nonumber \\
 & = & \sum_{i=0}^{s-1}\cos^{2}\left((2i+1)\frac{\pi}{4s+2}\right)\nonumber \\
 & = & \frac{1}{2}\sum_{i=0}^{s-1}\left(1+\cos\left(2(2i+1)\frac{\pi}{4s+2}\right)\right)\nonumber \\
 & = & \frac{s}{2}+\frac{\sin\left(\frac{s\pi}{2s+1}\right)\cos\left(\frac{s\pi}{2s+1}\right)}{2\sin\left(\frac{\pi}{2s+1}\right)}.\label{eq:sumOfxsquared}
\end{eqnarray}
Setting $\theta=\frac{\pi}{2s+1}$, we get $s=\frac{\pi}{2\theta}-\frac{1}{2}$
and then 
\[
\cos\left(\frac{s\pi}{2s+1}\right)=\sin\left(\frac{\theta}{2}\right)\text{ and }\sin\left(\frac{s\pi}{2s+1}\right)=\cos\left(\frac{\theta}{2}\right).
\]
Therefore, by means of trigonometric identities, we obtain $c=\pi/4\theta$.
Thus 
\begin{equation}
\omega^{2}=\frac{1}{c}=\frac{4\theta}{\pi}.\label{eq:normX}
\end{equation}

On the other hand, the Taylor series for $\left(\cos\left(\frac{(2r-1)\theta}{2}\right)-\cos\left(\frac{(2(r+k)-1)\theta}{2}\right)\right)^{2}$
is given by
\[
\frac{1}{4}\,k^{2}\left(2\,r+k-1\right)^{2}\theta^{4}+O\left(\theta^{6}\right).
\]
Therefore, 
\[
\left(x_{r}-x_{r+k}\right)^{2}=\frac{1}{\pi}\,k^{2}\left(2\,r+k-1\right)^{2}\theta^{5}+O\left(\theta^{7}\right),
\]
as required
\end{proof}
The last result enables us to prove the main Theorems of this section.
\begin{proof}
(Theorem \ref{thm:maxJump}) Recall that in our notation $x_{1}>x_{2}>\ldots>x_{s}$.
Fix $r=n^\alpha$ and $k=n^\beta$, and let $$R=\{(i,j):y_{j}<y_{i}\mbox{ and }i+k\leq j \mbox{ and } r\leq i\}.$$
Then 
\[
2n\Vert x-y\Vert^{2}=\sum_{i=1}^{n}\sum_{j=1}^{n}(x_{i}-y_{i})^{2}+(x_{j}-y_{j})^{2}\geq\sum_{(i,j)\in R}(x_{i}-y_{i})^{2}+(x_{j}-y_{j})^{2}.
\]

Since $y_{j}\geq y_{i}$ and $x_{j}<x_{i}$, the minimum contribution
that each term in the sum can provide happens when 
\[
y_{i}=y_{j}=\frac{x_{i}+x_{j}}{2}.
\]
Thus, we have 
\[
2n\left\Vert x-y\right\Vert ^{2}>\sum_{(i,j)\in R}\left(x_{i}-\frac{x_{i}+x_{j}}{2}\right)^{2}+\left(x_{j}-\frac{x_{i}+x_{j}}{2}\right)^{2}=\sum_{(i,j)\in R}\frac{\left(x_{i}-x_{j}\right)^{2}}{2}.
\]
Besides $i+k\geq j$, and then the last expression is minimum for
$j=i+k$. Then we can write
\[
2n\left\Vert x-y\right\Vert ^{2}>\sum_{(i,j)\in R}\frac{\left(x_{i}-x_{i+k}\right)^{2}}{2}.
\]
By Lemma \ref{lem:Xunitary}, for $n$ large enough 
\[
2n\left\Vert x-y\right\Vert ^{2}>\frac{1}{2\pi}\sum_{(i,j)\in R}\frac{k^{2}\left(2\,i+k-1\right)^{2}}{n^{5}}
>\frac{k^{2}}{2\pi n^{5}} \sum_{(i,j) R} i^2+k^{2}.
\]
By definition $r\leq i$. For $r=n^\alpha$ and $k=n^\beta$, we can write
$$2n\left\Vert x-y\right\Vert ^{2} > \frac{n^{2\beta}}{2\pi n^{5}} \sum_{(i,j) R} n^{2\alpha}+n^{2\beta}.  $$

%If we assume $|R|=\mathcal{O}(n^2)$, which means roughly all pairs are inverted, then we can write the lower bound
That gives the bound
$$2n\left\Vert x-y\right\Vert ^{2} > C_1 \frac{n^{2\beta}}{ n^{5}} |R|( n^{2\alpha}+n^{2\beta}),  $$
for an absolute constant $C_1>0$. Equivalently, we get
$$\left\Vert x-y\right\Vert ^{2} > C_1  |R| (n^{2(\alpha+\beta)-6}+n^{4\beta-6}).  $$

%For each set $R$ we will construct a set $R^*$ with the same cardinality of $R$ that satisfies
%$$\sum_{(i,j)\in R}\left(2\,i-1\right)^{2}>\sum_{(i,j)\in R^*}\left(2\,i-1\right)^{2}.$$
%That provides a universal lower bound on this sum in terms of the cardinality of $R$.
%To construct such set $R^*$ define
% $$R^*=\{(1,1+k),(2,2+k),\ldots,(|R|,|R|+k)\}.$$
%That gives us the bound
%\[2n\left\Vert x-y\right\Vert ^{2}>\frac{k^{4}}{2\pi n^{5}} \sum_{i=1}^{|R|}\left(2\,i-1\right)^{2}=\frac{k^{4}}{2\pi n^{5}} \frac{|R|(2|R|-1)(2|R|+1)}{3}.\]

%And for $|R|$  large enough we can find a absolute constant $C_2>0$ that gives the lower bound
%\[2n\left\Vert x-y\right\Vert ^{2}>C_2\frac{k^{4}}{ n^{5}} |R|^{3}.\] 
%Now let $k=n^{\beta}$. Then the lower bound becomes
%\[\left\Vert x-y\right\Vert ^{2}>C_2  |R|^{3}\frac{n^{4\beta}}{n^{6}}.\]

On the other hand, by Theorem \ref{thm:upperBound} there is a constant
$C_{0}>0$ such that with probability $1-n^{-3}$
\[
\left\Vert x-y\right\Vert ^{2}\leq C_{0}n^{-1}.
\]
 Combining these two inequalities, we obtain a constant $C_{2}>0$
such that 
$$|R|<C_2(n^{5-2(\alpha+\beta)}+n^{5-4\beta}).  $$
with probability $1-n^{-3}$.
That finishes the proof.
\end{proof}
\bigskip

\begin{proof}
(Theorem \ref{thm:maxJump2}) Recall that in our notation $x_{1}>x_{2}>\ldots>x_{s}$. 
Fix $k=n^\beta$, and let $$R=D_{k,1}=\{(i,j):y_{j}<y_{i}\mbox{ and }i+k\leq j  \}.$$
Then 
\[
2n\Vert x-y\Vert^{2}=\sum_{i=1}^{n}\sum_{j=1}^{n}(x_{i}-y_{i})^{2}+(x_{j}-y_{j})^{2}\geq\sum_{(i,j)\in R}(x_{i}-y_{i})^{2}+(x_{j}-y_{j})^{2}.
\]

Since $y_{j}\geq y_{i}$ and $x_{j}<x_{i}$, the minimum contribution
that each term in the sum can provide happens when 
\[
y_{i}=y_{j}=\frac{x_{i}+x_{j}}{2}.
\]
Thus, we have 
\begin{equation}
\label{eq:eqDiffSquared}
2n\left\Vert x-y\right\Vert ^{2}>\sum_{(i,j)\in R}\left(x_{i}-\frac{x_{i}+x_{j}}{2}\right)^{2}+\left(x_{j}-\frac{x_{i}+x_{j}}{2}\right)^{2}=\sum_{(i,j)\in R}\frac{\left(x_{i}-x_{j}\right)^{2}}{2}.
\end{equation}

If there are $n_i$ pairs $(i,j)$, then we can index $j$'s as
$(i, j_{i_{1}}), \dots (i,j_{i_{n_i}}), \text{ for } i=1,\ldots,n.$
Here, the total number of inverted pairs is $\sum_{i=1}^n n_i=|R|$.
Furthermore, for a fixed $i$ the we obtain the minimum $(x_i-x_j)^2$ whenever $j-i$ is minimum. Therefore, since $i+k\leq j$ for each $i$ the minimum sum
$$\sum_{t=1}^{n_i} \left(x_{i}-x_{j_{i_t}}\right)^{2}$$
occurs whenever $j_{i_{1}}=i+k, j_{i_{2}}=i+k+1, \dots, j_{i_{n_i}}=i+k+n_i-1$. Therefore, inequality (\ref{eq:eqDiffSquared}) becomes
\begin{equation*}
2n\left\Vert x-y\right\Vert ^{2} > \sum_{i=1}^n \sum_{t=0}^{n_i-1} \frac{\left(x_{i}-x_{i+k+t}\right)^{2}}{2}.
\end{equation*}
By Lemma \ref{lem:Xunitary}, for $n$ large enough 
\[
2n\left\Vert x-y\right\Vert ^{2}>\frac{1}{2\pi}\sum_{i=1}^n \sum_{t=0}^{n_i-1}   \frac{(k+t)^2(i+k+t)^2}{n^{5}}
>\frac{k^{2}}{2\pi n^{5}} \sum_{i=1}^n \sum_{t=0}^{n_i-1}  t^2.
\]
Therefore, for $n$ large enough there is a constant $C_1>0$ such that
\begin{equation}\label{eq:ncubed}
2n\left\Vert x-y\right\Vert ^{2} > C_1 k^{2} \sum_{i=1}^n n_i^3.
\end{equation}
Recall that two $p$-norms are related by $\left\Vert x\right\Vert_p\leq n^{\frac{1}{p}-\frac{1}{q}} \left\Vert x\right\Vert_q $.
Taking $p=1$ and $q=3$, we obtain that
$$ \sum_{i=1}^n n_i \leq n^{2/3} \sum_{i=1}^n n_i^3. $$
The last inequality together with $\sum_{i=1}^n n_i=|R|$ allows us to write inequality (\ref{eq:ncubed}) as
\begin{equation*}
2n\left\Vert x-y\right\Vert ^{2} > C_1 k^{2} |R|^3 n^{-2}.
\end{equation*}

On the other hand, by Theorem \ref{thm:upperBound} there is a constant
$C_{0}>0$ such that with probability $1-n^{-3}$
\[
\left\Vert x-y\right\Vert ^{2}\leq C_{0}n^{-1}.
\]
 Combining these two inequalities, we obtain a constant $C_{2}>0$
such that 
$$|R|<C_2 n^{\frac{7-2\beta}{3}}.  $$
with probability $1-n^{-3}$.
That finishes the proof.
\end{proof}
\bigskip

\bigskip
\begin{proof}
(Theorem \ref{thm:footrule}) Fix $k=k(n)=n^{4/5}$.  Define 
\[
R=\{(i,j):y_{j}<y_{i}\mbox{ and }i+k\leq j\},
\]
and
\[
R^{c}=\{(i,j):y_{j}<y_{i}\mbox{ and } j<i+k\}.
\]

By Theorem \ref{thm:maxJump2}, taking $\beta=4/5$, there exists a
constant $C$ so that, for large enough $n$, 
\[
|R|=D_{k}(y)\leq Cn^{9/5}.
\]
On the other hand, for each index $i$ there can be at most $k$ pairs
$(i,j)$ that occur in $R^{c}$, so 
\[
|R^{c}|\leq kn=n^{9/5}.
\]
 Therefore, $D(\sigma)=|R|+|R^{c}|\leq(C+1)n^{9/5}$. 
\end{proof}
We finish this section showing that there is a vector that gives $D(\sigma)=\mathcal{O}(n^{8/5})$.
\begin{thm}
Let $x$ be the unitary eigenvector for $\lambda_{2}\left(M\right)$,
and let $\sigma$ be the permutation derived from any vector $\hat{x}$,
such that $\hat{x}_{\sigma(1)}\geq\hat{x}_{\sigma(2)}\geq\ldots\geq\hat{x}_{\sigma(n)}$.
Then, there is a vector $\hat{x}$ with a permutation satisfying $D(\sigma)=\mathcal{O}(n^{8/5})$
and $\left\Vert \hat{x}-x\right\Vert \leq\mathcal{O}(n^{-1/2})$.\end{thm}
\begin{proof}
We prove it by constructing such vector. Denote the standard basis
vectors of $\mathbb{R}^{n}$ by $e_{i}$. Let $1\leq k<n/2$ be an
integer and define the vector $P=1/\sqrt{k}\left[\begin{array}{cccccc}
1 & \ldots & 1 & 0 & \ldots & 0\end{array}\right]^{T}$, where the entries $i=1,\ldots,k$ contain $1$ and $0$ elsewhere. 

We call $x*$ the projection of $x$ onto the subspace $S$ with orthonormal
basis given by 
\[
\left\{ P,e_{i}\text{ for }i=k+1,\ldots,n\right\} .
\]

Notice that, $\left\Vert x*\right\Vert <\left\Vert x\right\Vert =1$.
Define $y*$ as the vector obtained from  $x*$ by rescaling it in
order to get $\left\Vert y*\right\Vert =1$. The vector $y*$ will
be the required $\hat{x}$ of the Theorem, for some $k$ to be chosen
later.

Then, since $x*$ is an orthogonal projection onto $S$, we have 

\begin{eqnarray}
\left\Vert x-y*\right\Vert ^{2} & = & \left\Vert x-x*\right\Vert ^{2}+\left\Vert x*-y*\right\Vert ^{2}\nonumber \\
 & = & \left\Vert x-x*\right\Vert ^{2}+\left(1-\left\Vert x*\right\Vert \right)^{2}.\label{eq:pythagoras}
\end{eqnarray}

Also, the orthogonal projection matrix of a vector $x$ onto $S$
provides us the relation $x*=UU^{T}x$, where $U=\left[\begin{array}{cccc}
P & e_{k+1} & \ldots & e_{n}\end{array}\right]$. Besides, 
\[
PP^{T}x=1/k\left[\begin{array}{cccccc}
b & \ldots & b & 0 & \ldots & 0\end{array}\right]^{T},
\]
where 
\[
b=\sum_{i=1}^{k}x_{i}.
\]
Thus, we can write
\begin{eqnarray*}
x-x* & = & \left(I-UU^{T}\right)x\\
 & = & \left[\begin{array}{cccccc}
x_{1}-\frac{b}{k} & \ldots & x_{k}-\frac{b}{k} & 0 & \ldots & 0\end{array}\right]^{T}.
\end{eqnarray*}
That shows 

\begin{eqnarray}
\left\Vert x-x*\right\Vert ^{2} & = & \sum_{i=1}^{k}\left(\frac{b}{k}-x_{i}\right)^{2}\nonumber \\
 & = & \sum_{i=1}^{k}\frac{b^{2}}{k^{2}}-2\frac{b}{k}x_{i}+x_{i}^{2}\nonumber \\
 & = & \frac{b^{2}}{k}-2\frac{b}{k}\sum_{i=1}^{k}x_{i}+\sum_{i=1}^{k}x_{i}^{2}\nonumber \\
 & = & \sum_{i=1}^{k}x_{i}^{2}-\frac{b^{2}}{k}\nonumber \\
 & = & \sum_{i=1}^{k}x_{i}^{2}-\frac{\left(\sum_{i=1}^{k}x_{i}\right)^{2}}{k}.\label{eq:expNormofX}
\end{eqnarray}

Furthermore, Theorem \ref{thm:eigA} provides the expressions we need
to compute $\left\Vert x-x*\right\Vert ^{2}$: we have $x_{j}=\omega\cos\left(\frac{(2j-1)\pi}{4s+2}\right)$,
for $j=1,\ldots,s$, where $\omega$ is a constant such that $\left\Vert x\right\Vert =1$. 

Then, we can write

\begin{eqnarray*}
b & = & \sum_{i=1}^{k}\omega\cos\left(\frac{(2i-1)\pi}{4s+2}\right)\\
 & = & \sum_{i=0}^{k-1}\omega\cos\left((2r+2i+1)\frac{\pi}{4s+2}\right).
\end{eqnarray*}
Additionally, we make use of the trigonometric identities 
\[
\sum_{i=0}^{k-1}\sin(2\alpha i+\beta)=\frac{\sin\left(k\alpha\right)\sin\left(\beta+\left(k-1\right)\alpha\right)}{\sin\alpha}\text{ and }
\]
\[
\sum_{i=0}^{k-1}\cos(2\alpha i+\beta)=\frac{\sin\left(k\alpha\right)\cos\left(\beta+\left(k-1\right)\alpha\right)}{\sin\alpha}.
\]
Then, $b$ can be rewritten as 
\begin{equation}
b=\omega\frac{\sin\left(\frac{k\pi}{4s+2}\right)\cos\left(\frac{2r\pi}{4s+2}+\frac{k\pi}{4s+2}\right)}{\sin\left(\frac{\pi}{4s+2}\right)}.\label{eq:sumOfx}
\end{equation}
Also, defining $c=\sum_{i=1}^{k}x_{i}^{2}$, we obtain 
\begin{eqnarray}
c & = & \sum_{i=0}^{k-1}\left(\omega\cos\right)^{2}\left((2i+1)\frac{\pi}{4s+2}\right)\nonumber \\
 & = & \frac{\omega^{2}}{2}\sum_{i=0}^{k-1}\left(1+\cos\left(2(2i+1)\frac{\pi}{4s+2}\right)\right)\nonumber \\
 & = & \omega^{2}\frac{k}{2}+\omega^{2}\frac{\sin\left(\frac{2k\pi}{4s+2}\right)\cos\left(\frac{2k\pi}{4s+2}\right)}{2\sin\left(\frac{2\pi}{4s+2}\right)}.\label{eq:sumOfxsquared-1}
\end{eqnarray}

Thus, to get an expression for $\left\Vert x-x*\right\Vert ^{2}$
we put together equations (\ref{eq:expNormofX}),(\ref{eq:sumOfx}),
and (\ref{eq:sumOfxsquared-1}) and we obtain 
\begin{eqnarray}
\left\Vert x-x*\right\Vert ^{2} & = & c-\frac{b^{2}}{k}\nonumber \\
 & = & \omega^{2}\frac{k}{2}+\omega^{2}\frac{\sin\left(\frac{2k\pi}{4s+2}\right)\cos\left(\frac{2k\pi}{4s+2}\right)}{2\sin\left(\frac{2\pi}{4s+2}\right)}\nonumber \\
 &  & -\frac{\omega^{2}}{k}\left(\frac{\sin\left(\frac{k\pi}{4s+2}\right)\cos\left(\frac{k\pi}{4s+2}\right)}{\sin\left(\frac{\pi}{4s+2}\right)}\right)^{2}.\label{eq:x-xx}
\end{eqnarray}

Now, in view of equation (\ref{eq:pythagoras}), to get an expression
for $\left\Vert x-y*\right\Vert $ we need an expression for $\left\Vert x*\right\Vert $.
But,
\begin{eqnarray*}
\left\Vert x*\right\Vert ^{2} & = & \left\Vert x-(x-x*)\right\Vert ^{2}\\
 & = & \left\Vert x\right\Vert ^{2}-2\left\langle x,x-x*\right\rangle +\left\Vert x-x*\right\Vert ^{2}.
\end{eqnarray*}
Furthermore, we have
\begin{eqnarray*}
\left\langle x,x-x*\right\rangle  & = & \sum_{i=1}^{k}x_{i}\left(x_{i}-\frac{b}{k}\right)\\
 & = & \sum_{i=1}^{k}x_{i}^{2}-\frac{b}{k}\sum_{i=1}^{k}x_{i}\\
 & = & c-\frac{b^{2}}{k}.
\end{eqnarray*}
Therefore, we can write
\begin{eqnarray*}
\left\Vert x*\right\Vert ^{2} & = & \left\Vert x\right\Vert ^{2}-2\left(c-\frac{b^{2}}{k}\right)+\left\Vert x-x*\right\Vert ^{2}.
\end{eqnarray*}
And by equation (\ref{eq:x-xx}), we get
\begin{eqnarray*}
\left\Vert x*\right\Vert ^{2} & = & \left\Vert x\right\Vert ^{2}-2\left(c-\frac{b^{2}}{k}\right)+c-\frac{b^{2}}{k}\\
 & = & \left\Vert x\right\Vert ^{2}-c+\frac{b^{2}}{k}\\
 & = & \left\Vert x\right\Vert ^{2}-\left\Vert x-x*\right\Vert ^{2}.
\end{eqnarray*}
Now, we plug this expression in $\left(1-\left\Vert x*\right\Vert \right)^{2}$
and use the fact that $\left\Vert x\right\Vert =1$ to obtain 

\begin{eqnarray}
\left(1-\left\Vert x*\right\Vert \right)^{2} & = & 1-2\left\Vert x*\right\Vert +\left\Vert x*\right\Vert ^{2}\nonumber \\
 & = & 1-2\sqrt{\left\Vert x\right\Vert ^{2}-\left\Vert x-x*\right\Vert ^{2}}+\left\Vert x\right\Vert ^{2}-\left\Vert x-x*\right\Vert ^{2}\\
 & = & 2-2\sqrt{1-\left\Vert x-x*\right\Vert ^{2}}-\left\Vert x-x*\right\Vert ^{2}\\
 & = & 2-2\sqrt{1-\left(c-\frac{b^{2}}{k}\right)}-\left\Vert x-x*\right\Vert ^{2}.\label{eq:1-xx}
\end{eqnarray}

Finally, equations (\ref{eq:pythagoras}), (\ref{eq:x-xx}), and (\ref{eq:1-xx}),
give us the expression
\begin{eqnarray}
\left\Vert x-y*\right\Vert ^{2} & = & \left\Vert x-x*\right\Vert ^{2}+\left(1-\left\Vert x*\right\Vert \right)^{2}\nonumber \\
 & = & 2-2\sqrt{1-\left(c-\frac{b^{2}}{k}\right)}.\label{eq:minc-b}
\end{eqnarray}
Thus, $\left\Vert x-y*\right\Vert ^{2}$ is completely determined
by the function prescribed in equation (\ref{eq:x-xx}). Therefore,
we get the function 
\begin{eqnarray*}
f(k,s) & = & \frac{k}{2}+\frac{\sin\left(\frac{2k\pi}{4s+2}\right)\cos\left(\frac{2k\pi}{4s+2}\right)}{2\sin\left(\frac{2\pi}{4s+2}\right)}\\
 &  & -\frac{1}{k}\left(\frac{\sin\left(\frac{k\pi}{4s+2}\right)\cos\left(\frac{k\pi}{4s+2}\right)}{\sin\left(\frac{\pi}{4s+2}\right)}\right)^{2}.
\end{eqnarray*}
Therefore, we obtain 
\begin{equation}
\left\Vert x-y*\right\Vert ^{2}=2-2\sqrt{1-\omega^{2}f(k,s)},\label{eq:normnorm}
\end{equation}
where $\omega=1/\left\Vert x\right\Vert $. 

Let $\theta=\frac{\pi}{4s+2}$. We look to the Taylor series of $f$
at $\theta=0$. With aid of a computer algebra system, we obtain
\begin{eqnarray*}
f(k,s) & = & \theta^{4}\left(\vphantom{\frac{}{}}\right.{\frac{16}{45}}\,{k}^{5}-\frac{4}{9}\,{k}^{3}+{\frac{4}{45}}\,k\left.\frac{}{}\right)+\mathcal{O}\left(\theta^{6}\right).
\end{eqnarray*}
For the value $\omega^{2}$ we use the same estimative as in Lemma
\ref{lem:Xunitary} where we obtained equation (\ref{eq:normX}).
Thus equation (\ref{eq:normnorm}) is 
\[
\left\Vert x-y*\right\Vert ^{2}=\mathcal{O}\left(\theta^{5}k^{5}\right)=\mathcal{O}\left(\frac{k^{5}}{n^{5}}\right).
\]

Notice that the vector $y*$ has its first $k$ entries in the wrong
order, and the remaining equals to the corresponding entries in $x$,
so they are correct. Now if we choose $k=n^{4/5}$ and take $\widehat{x}=y*$,
the permutation $\sigma$ in the Theorem is such that 
\[
D(\sigma)=\binom{k}{2}=\mathcal{O}(n^{8/5}).
\]

Besides, we have 
\[
\left\Vert x-\widehat{x}\right\Vert ^{2}=\mathcal{O}\left(\frac{(n^{4/5})^{5}}{n^{5}}\right)\leq\mathcal{O}\left(\frac{1}{n}\right).
\]
That finishes the proof.
\end{proof}

\section{Final remarks}

We finish this paper by highlighting the generality of our method.
Once the eigenvectors that reveal the structure of the model graph
are identified, we can use a similar technique to recover the structure
of the random graph with high probability.

First, the distance between the eigenvector of the random matrix and
the model matrix can be always bounded by general results from random
matrix theory. In our method we can significantly refine these bounds
by using the eigenvalues. Second, to quantify the vertices that are
incorrectly placed, we rely on series expansions for the entries of
the eigenvector (see Section \ref{sec:qualitative}). It is worth
mentioning that the trade off between how far vertices can jump and
the total amount, i.e., the proportion of vertices that are shifted
significantly from their true positions, is negligible. This behaviour
seems to be a general feature for these kind of problems, as the proof
of Theorem \ref{thm:maxJump} reveals.

Finally, this paper serves as proof of concept, and in a forthcoming
work we will deal with different geometric models, such as grids,
rings, toroids, product of graphs, etc. At this point, it is clear
one need to find analytic expressions for the eigenvectors of interest
in order to reveal the structure of a random graph.

\end{document}